\documentclass[a4paper,11pt,reqno]{amsart}

\usepackage[margin=20mm,includehead,includefoot]{geometry}
\usepackage[utf8]{inputenc}
\usepackage[T1]{fontenc}
\usepackage{lmodern,microtype}
\usepackage{upref}
\usepackage{mathtools,amssymb,amsthm,amsmath}
\usepackage[foot]{amsaddr}
\usepackage{enumitem}
\usepackage[textsize=footnotesize,disable]{todonotes}
\usepackage{graphicx}
\usepackage{hyperref}
\usepackage{wrapfig}
\usepackage[margin=1cm]{caption}
\usepackage{bookmark}
\usepackage{mycommands}

\graphicspath{{./graphics/}}

\renewenvironment{quote}{%
  \list{}{%
    \leftmargin4mm   % this is the adjusting screw
    \rightmargin\leftmargin
  }
  \item\relax
}
{\endlist}

\makeatletter
\let\old@setaddresses\@setaddresses
\def\@setaddresses{\bigskip{\parindent 0pt\let\scshape\relax\let\ttfamily\relax\old@setaddresses}}
\makeatother


\hypersetup{
  pdftitle={On Hamilton cycles in graphs defined by intersecting set systems},
  pdfauthor={Torsten M\"utze}
}

\title{On Hamilton cycles in graphs \\ defined by intersecting set systems}

\author{Torsten M\"utze}
\address[Torsten M\"utze]{Department of Computer Science, University of Warwick, United Kingdom \& Department of Theoretical Computer Science and Mathematical Logic, Charles University, Prague, Czech Republic}
\email{torsten.mutze@warwick.ac.uk}

\thanks{This work was supported by Czech Science Foundation grant GA~22-15272S}

\begin{document}

\maketitle

\section{Introduction}
\label{sec:intro}

\begin{wrapfigure}{r}{0.45\textwidth}
\includegraphics[scale=0.25]{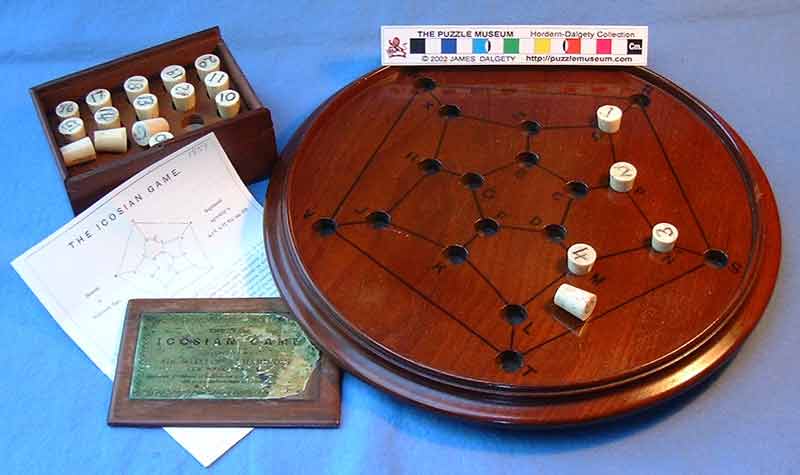}
\caption{Hamilton's Icosian game.}
\vspace{-2mm}
\end{wrapfigure}
In 1857 the Irish mathematician William Rowan Hamilton invented a puzzle whose goal is to find a cycle in the graph of the dodecahedron that visits every vertex exactly once.
He dubbed it the `Icosian game', as the resulting cycle has exactly twenty (`icosa' in ancient Greek) edges and vertices.
In honor of Hamilton, a cycle that visits every vertex of a graph exactly once is now called a \emph{Hamilton cycle}. 
The dodecahedron has the interesting property that it looks the same from the point of view of any vertex.
Formally, it is \emph{vertex-transitive}, i.e., any two vertices can be mapped onto each other by an automorphism of the graph.
In 1970 Lov\'asz raised a conjecture which can be considered a highly advanced version of the Icosian game.
Specifically, he conjectured that every connected vertex-transitive graph admits a Hamilton cycle, apart from five exceptional graphs, which are a single edge, the Petersen graph, the Coxeter graph, and the graphs obtained from the latter two by replacing every vertex by a triangle.

\begin{wrapfigure}{r}{0.45\textwidth}
\centering
\includegraphics{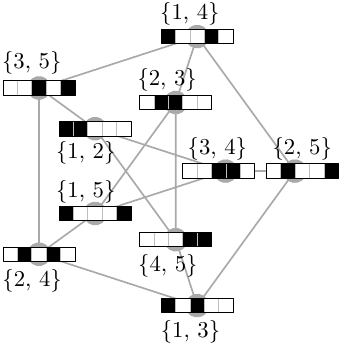}
\caption{The Petersen graph as intersection graph of all 2-element subsets of~$\{1,2,3,4,5\}$, with edges connecting disjoint sets.
In the corresponding bitstring representations, 0s are drawn as white squares and 1s as black squares.}
\label{fig:petersen}
\vspace{-4mm}
\end{wrapfigure}
The Petersen graph, shown in Figure~\ref{fig:petersen}, is a 10-vertex graph that serves as a famous example and counterexample for many problems in graph theory, and we will encounter it again later in this article.
A variant of Lov\'asz' conjecture asserts that every connected vertex-transitive graph admits a \emph{Hamilton path} (without exceptions), i.e., we still need to visit every vertex exactly once, but the first and last vertex of the tour need not be adjacent.
While the Icosian game is about one particular graph, Lov\'asz' conjecture talks about infinitely many graphs, which on the one hand are very constrained, but on the other hand very hard to get your hands on.
One can easily construct numerous explicit families of vertex-transitive graphs for which it is not known whether they are Hamiltonian.
An important example for this are \emph{Cayley graphs}, which are defined for a group and a set of generators as the graph that has as vertices all group elements, and whose edges arise by multiplication with a generator.
Even for particular groups like the symmetric group, and for generators with certain properties (like two generators, or generators that are involutions), only partial results are known.

In this expository article, we consider another rich family of vertex-transitive graphs, which features prominently throughout combinatorics, namely graphs defined by intersecting set systems.
We give an overview of the many beautiful techniques and ingredients devised during the past 40~years that establish the existence of Hamilton cycles in those graphs, thus settling interesting special cases of Lov\'asz' conjecture.
Our discussion starts with the easy instances and ends with the hardest and most general ones, and it emphasizes the key obstacles on this journey.

\subsection{The hypercube}
\label{sec:cube}

\begin{figure}[h!]
\begin{tabular}{cc}
\includegraphics{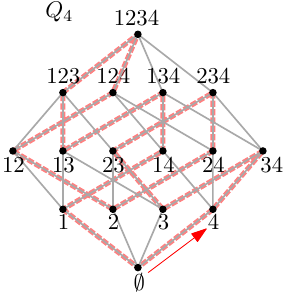} \hspace{5mm} &
\includegraphics{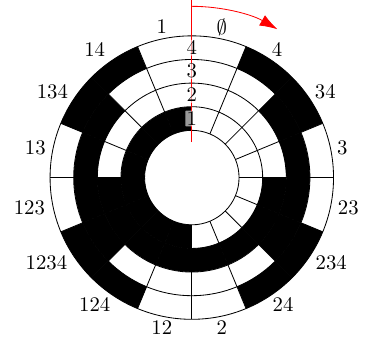}
\end{tabular}
\caption{Left: The 4-dimensional hypercube~$Q_4$ and one of its Hamilton cycles, the binary reflected Gray code (highlighted edges); Right: bitstring representation of the cycle ($0=\text{white}$ and $1=\text{black}$) with vertices in clockwise order starting at 12 o'clock.
When printing the sets, curly brackets and commas are omitted for simplicity.
}
\label{fig:q4}
\end{figure}

The starting point is the Boolean lattice, i.e., the inclusion order on all subsets of $[n]:=\{1,2,\ldots,n\}$; see Figure~\ref{fig:q4}.
The cover relations of this poset are pairs of sets $X,Y$ which differ in a single element, i.e., $Y=X\cup\{i\}$ for some $i\in[n]$, and the corresponding cover graph~$Q_n$ is the well-known \emph{$n$-dimensional hypercube}.
It is convenient to encode the vertices of~$Q_n$ as bitstrings of length~$n$, by considering the characteristic vector of each set, which has the $i$th bit equal to 1 if and only if the element~$i$ is contained in the set.
With this encoding, edges of~$Q_n$ connect exactly pairs of vertices that differ in a single bit.
This graph is vertex-transitive, and a Hamilton cycle can be found easily, by applying the following simple greedy rule discovered by Williams: Start at an arbitrary vertex, and repeatedly flip the rightmost bit that creates a previously unvisited vertex.
The resulting cycle is known as \emph{binary reflected Gray code}, named after Bell Labs researcher Frank Gray, and it is shown in Figure~\ref{fig:q4} for $n=4$.

\section{The middle levels conjecture}

The \emph{$k$th level} of~$Q_n$ is the set of all vertices with exactly~$k$ many~1s.
In terms of subsets in the Boolean lattice, these are the subsets of size exactly~$k$.
Now consider the hypercube of odd dimension~$2k+1$, and the subgraph induced by the middle two levels~$k$ and~$k+1$; see Figure~\ref{fig:m2}.
We denote this subgraph by~$M_k$, and we note that it is vertex-transitive.
Havel~\cite{MR737021} in~1983, and independently Buck and Wiedemann~\cite{MR737262} in~1984 conjectured that $M_k$ has a Hamilton cycle for all~$k\geq 1$, and this problem became known as \emph{middle levels conjecture}.
It turns out that this problem is considerably harder than finding a Hamilton cycle in the entire cube.
Thus, the conjecture is a prime example of an easy-to-state combinatorial proposition that one feels should be easy to prove at first sight, but that turns out to be surprisingly intricate upon further investigation.
Also, it is an explicit instance of Lov\'asz' conjecture, and failing to prove it for this particular family of graphs~$M_k$ shows how little we understand about the general problem.
Consequently, the middle levels conjecture has attracted a lot of attention in the literature (see e.g., \cite{MR962224,MR1350586,MR1329390}, and it is mentioned in several popular books, in particular in Diaconis and Graham's book `Magical mathematics', and in Winkler's book `Mathematical puzzles: a connoisseur's collection'.
Furthermore, Knuth gave the middle levels conjecture the highest difficulty rating (49/50) among all open problems in his book `The Art of Computer Programming Vol.~4A'.
In a recent survey, Gowers comments on the conjecture as follows:
\begin{quote}
`If one starts trying to build a Hamilton cycle in $M_k$, one runs into the problem of having too much choice and no obvious way of making it. (A natural thing to try to do is find some sort of inductive construction, but a lot of people have tried very hard to do this, with no success---a natural pattern just doesn't seem to emerge after the first few small cases.)'
\end{quote}

\subsection{Kierstead-Trotter matchings}

A \emph{matching} in a graph is a set of pairwise disjoint edges, and a matching is \emph{perfect} if it includes every vertex.
Kierstead and Trotter~\cite{MR962224} in~1988 suggested to tackle the middle levels conjecture by taking the union of two edge-disjoint perfect matchings in~$M_k$, with the hope that their union forms the desired Hamilton cycle.

\begin{wrapfigure}{r}{0.48\textwidth}
\begin{tabular}{c}
\hspace{-3mm}\includegraphics[page=1]{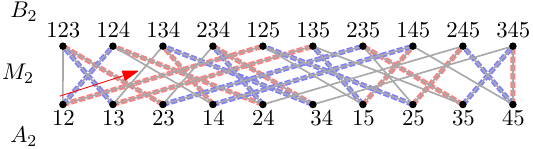} \\
\hspace{10mm}\includegraphics{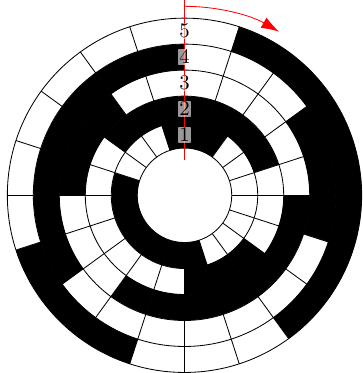}
\end{tabular}
\caption{Top: The middle levels graph~$M_2$ and the perfect matchings~$X_f$ and~$X_g$ (red and blue, respectively) whose union is a Hamilton cycle; Bottom: bitstring representation of the cycle.}
\label{fig:m2}
\end{wrapfigure}
The two matchings they considered can be described explicitly as follows:
We write $A_k$ and~$B_k$ for the vertices in levels~$k$ and~$k+1$ of~$Q_{2k+1}$, respectively, i.e., these are the two partition classes of~$M_k$.
A \emph{Dyck word} is a bitstring with the same number of~1s and~0s and the property that every prefix contains at least as many~1s as~0s.
We write $D_k$ for the set of all Dyck words of length~$2k$.
Furthermore, for any string~$x$ and any integer~$i$, we write $\sigma^i(x)$ for the cyclic left shift of~$x$ by $i$ steps.
Note that every vertex~$x\in A_k$ can be written uniquely as $x=\sigma^i(y\,0)$ for some $y\in D_k$ and $0\leq i<2k+1$.
Indeed, the counting works out correctly, as the number of Dyck words is the Catalan number~$C_k=|D_k|=\frac{1}{k+1}\binom{2k}{k}$, the number of cyclic shifts is~$2k+1$, and $(2k+1)\cdot C_k=\binom{2k+1}{k}=|A_k|$.
We then define $f(x):=\sigma^i(y\,1)\in B_k$.
Furthermore, we decompose~$y$ uniquely as $y=1\,u\,0\,v$ for $uv\in D_{k-1}$, and we define $g(x):=\sigma^i(1\,u\,1\,v\,0)\in B_k$.
It is easy to show that~$f$ and~$g$ are bijections between~$A_k$ and~$B_k$, and therefore $X_f:=\{(x,f(x))\mid x\in A_k\}$ and $X_g:=\{(x,g(x))\mid x\in A_k\}$ are two edge-disjoint perfect matchings between the two partition classes~$A_k$ and~$B_k$ of~$M_k$.
The union of~$X_f$ and~$X_g$ thus forms a \emph{cycle factor}~$F_k:=X_f\cup X_g$ in the graph~$M_k$, i.e., a collection of disjoint cycles that together visit all vertices.
In fact, for $k=2$ we are lucky and~$F_k$ is a single cycle, i.e., a Hamilton cycle in~$M_k$; see Figure~\ref{fig:m2}.
Unfortunately, our luck ends for larger values of~$k$, as the number of cycles of~$F_k$ for $k=1,2,\ldots$ is $1,1,2,3,6,14,34,95,280,854,\ldots$.

\subsection{Interpretation of the cycles as plane trees}

The middle levels conjecture was finally solved by~M\"utze~\cite{MR3483129} in 2016, who later provided a 2-page proof~\cite{muetze:23}.
The short proof picks up on Kierstead and Trotter's argument as follows:
The counting sequence for the number of cycles of~$F_k$ is OEIS sequence~A002995, which counts plane trees with $k$ edges, hinting at a bijection between plane trees with $k$ edges and cycles in~$F_k$; see Figure~\ref{fig:aux}.

\begin{wrapfigure}{r}{0.6\textwidth}
\vspace{-5mm}
\includegraphics[page=2]{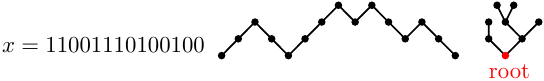}
\caption{Bijection between Dyck words (left), Dyck paths (middle) and ordered rooted trees (right).
Our trees have the root (red) at the bottom and grow upwards.}
\vspace{-2mm}
\label{fig:bij}
\end{wrapfigure}
Indeed, this bijection follows from the definition of~$f$ and~$g$, and it also uses the following well-known bijection between Dyck paths with $2k$ steps and ordered rooted trees with $k$ edges; see Figure~\ref{fig:bij}:
Given a Dyck word~$x$, we consider the corresponding Dyck path, which has an $\ustep$-step for every 1-bit and a $\dstep$-step for every 0-bit of~$x$.
We then squeeze this Dyck path together, gluing every pair of an $\ustep$-step and $\dstep$-step on the same height that `see' each other (i.e., that have no other step at the same height in between them) together to form an edge of the ordered rooted tree.

\begin{wrapfigure}{r}{0.4\textwidth}
\includegraphics[page=3]{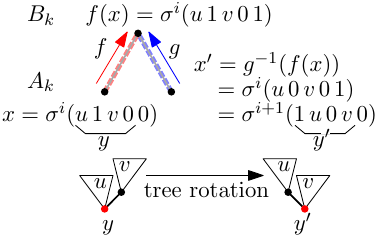}
\caption{Two steps along a cycle of~$F_k$ correspond to a tree rotation and cyclic left shift.}
\end{wrapfigure}
For a vertex~$x\in A_k$, we consider the vertex~$x'\in A_k$ that is two steps away along a cycle of~$F_k$, i.e., $x'$ is obtained by traversing one edge of~$X_f$ and one edge of~$X_g$, so $x'=g^{-1}(f(x))$.
Specifically, for $x=\sigma^i(u\,1\,v\,0\,0)\in A_k$ with $uv\in D_{k-1}$ and $0\leq i<2k+1$ we have $x':=g^{-1}(f(x))=\sigma^i(u\,0\,v\,0\,1)=\sigma^{i+1}(1\,u\,0\,v\,0)$.
We now define $y:=u\,1\,v\,0\in D_k$ and $y':=1\,u\,0\,v\in D_k$, and we consider the ordered rooted trees corresponding to the Dyck words~$y$ and~$y'$, which differ in a tree rotation.
Consequently, walking two steps along a cycle of~$F_k$ corresponds to a tree rotation and a cyclic left shift of the bitstring.
As $2k+1$ (the number of shifts) and $k$ (the number of tree edges) are coprime, all cyclic shifts of every tree lie on the same cycle.
As the equivalence classes of ordered rooted trees under tree rotation are precisely plane trees, obtained by `forgetting' the root, we obtain the desired bijection.
Now that we have a nice combinatorial interpretation of the cycles in the factor~$F_k$, all that is left to do is to join the cycles to a single Hamilton cycle.

\subsection{Gluing cycles}

\begin{wrapfigure}{r}{0.42\textwidth}
\vspace{-2mm}
\includegraphics{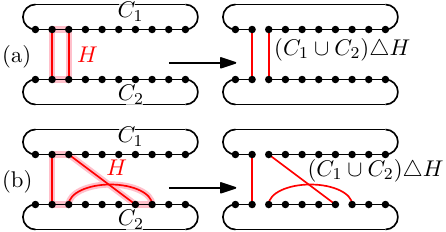}
\caption{Gluing cycles (red) join cycles from the cycle factor (black).}
\label{fig:gluing}
\end{wrapfigure}
The proof proceeds by gluing the cycles of the factor~$F_k$ together via small gluing cycles.
Given a cycle factor in a graph, a \emph{gluing cycle}~$H$ for a set $C_1,\ldots,C_\ell$ of cycles from the factor has every second edge in common with one of the cycles~$C_i$, in such a way that the symmetric difference of the edge sets of $(C_1\cup\cdots\cup C_\ell)\triangle H$ is a single cycle on the same vertex set as $C_1,\ldots,C_\ell$.
In the simplest case $\ell=2$ the cycle~$H$ is a 4-cycle that joins two cycles~$C_1,C_2$ as in Figure~\ref{fig:gluing}~(a).
Unfortunately, the middle levels graph~$M_k$ has no 4-cycles at all.
However, we can use a 6-cycle $H$ that intersects $C_1,C_2$ as shown in Figure~\ref{fig:gluing}~(b), and obtain the same effect.

\begin{wrapfigure}{r}{0.42\textwidth}
\includegraphics{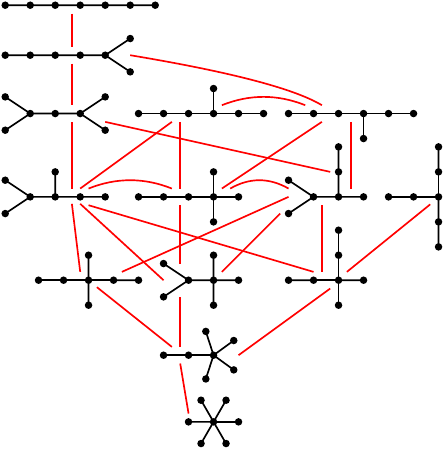}
\caption{Auxiliary graph~$\cH_k$ on plane trees with $k=6$ edges.}
\vspace{-6mm}
\label{fig:aux}
\end{wrapfigure}

As it turns out, there is a large collection of such gluing 6-cycles for the factor~$F_k$, and all these gluing cycles are edge-disjoint, i.e., the gluing operations do not interfere with each other.
Furthermore, each gluing cycle that joins cycles $C_1$ and~$C_2$ corresponds to a local modification of the plane trees associated with~$C_1$ and~$C_2$, which consists in removing a leaf from the plane tree, and reattaching it to a neighboring vertex.
We can thus define an auxiliary graph $\cH_k$ which has as nodes all plane trees with $k$ edges, corresponding to the cycles of the factor~$F_k$, and which has edges between pairs of plane trees that differ in such a local modification, corresponding to gluing cycles; see Figure~\ref{fig:aux}.
It remains to show that the auxiliary graph~$\cH_k$ is connected, which is done by showing that every plane tree can be transformed into a star by a sequence of local modifications as described before.

We thus reduced the problem of proving that $M_k$ has a Hamilton cycle to proving that the auxiliary graph~$\cH_k$ is connected, which is much easier, as nodes and edges of~$\cH_k$ have a nice combinatorial interpretation.
This two-step approach of building a Hamilton cycle (cycle factor+gluing) and the corresponding reduction to a spanning tree problem is very powerful, and has been employed also in several of the proofs discussed later.

\section{Bipartite Kneser graphs}

For integers $k\geq 1$ and $n\geq 2k+1$, the \emph{bipartite Kneser graph~$H_{n,k}$} has as vertices all $k$-element and $(n-k)$-element subsets of~$[n]$, and an edge between any two subsets~$X$ and~$Y$ with~$X\seq Y$.
The bipartite Kneser graph~$H_{n,k}$ is the cover graph of the subposet of the Boolean lattice of~$[n]$ induced by the levels~$k$ and~$n-k$.
In particular, $H_{2k+1,k}=M_k$ is the middle levels graph, i.e., bipartite Kneser graphs generalize the middle levels graphs.
Furthermore, bipartite Kneser graphs are vertex-transitive, which makes them interesting test cases for Lov\'asz' conjecture.
Simpson~\cite{MR1152123} and independently Roth conjectured in 1991 that all bipartite Kneser graphs admit a Hamilton cycle.
Note that the degrees of~$H_{n,k}$ are large when $n$ is large w.r.t.~$k$, i.e., intuitively, the middle levels case $n=2k+1$ is the sparsest and hardest one, whereas the denser cases $n>2k+1$ should be easier to prove.
The densest graph $H_{n,1}$ is the cover graph of what poset theorists call the `standard example', namely a complete bipartite graph minus a perfect matching.
Indeed, there has been considerable work on establishing that sufficiently dense bipartite Kneser graphs~$H_{n,k}$ have a Hamilton cycle.
Once the sparsest case $n=2k+1$ was established with the proof of the middle levels conjecture, the Hamiltonicity of all~$H_{n,k}$ was shown shortly thereafter by M\"utze and Su~\cite{MR3759914} in~2017.
In fact, their proof is a 5-page inductive argument, which uses the sparsest case~$n=2k+1$ as a basis.

\subsection{Havel's construction and its subsequent refinement}

\begin{wrapfigure}{r}{0.4\textwidth}
\vspace{-5mm}
\includegraphics{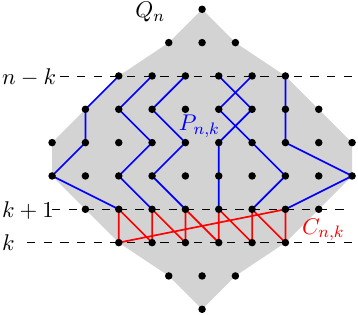}
\caption{Structures in~$Q_n$ used for the proof that $H_{n,k}$ has a Hamilton cycle.}
\label{fig:bkneser}
\end{wrapfigure}
We write~$Q_{n,k}$ for the subgraph of the hypercube~$Q_n$ induced by levels~$k$ and~$k+1$.
Already Havel~\cite{MR737021} in his 1983 paper considered the following strengthening of the middle levels conjecture:
For any $k\geq 1$ and $n\geq 2k+1$, there is a cycle~$C_{n,k}$ in~$Q_{n,k}$ that visits all vertices in level~$k$, i.e., in the smaller of the two partition classes, shown in red in Figure~\ref{fig:bkneser}.
For $n=2k+1$ both partition classes have the same size and $C_{n,k}$ is a Hamilton cycle in~$M_k$, i.e., this statement is the middle levels conjecture.
For $n>2k+1$ Havel proved this statement by an easy induction:
Indeed, we can split~$Q_{n,k}$ into two subgraphs~$Q_{n-1,k}$ and~$Q_{n-1,k-1}$, by partitioning all vertices according to the value of the last bit.
Using induction, we can glue together the two cycles~$C_{n-1,k}$ and~$C_{n-1,k-1}$ to obtain~$C_{n,k}$.
Note that this approach fails for $n=2k+1$, as in this case~$Q_{n-1,k}$ lies above the middle and~$Q_{n-1,k-1}$ lies below the middle, so the size difference between lower or upper level is opposite in both parts.

\begin{wrapfigure}{r}{0.42\textwidth}
\vspace{-8mm}
\centerline{
\setlength{\tabcolsep}{2pt}
\renewcommand{\arraystretch}{1.5}
\begin{tabular}{cccccccc}
(a1) & (a2) & (a3) & (a4) & (b1) & (b2) & (b3) & (b4) \\[-5mm]
\raisebox{-\height}{\includegraphics{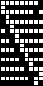}} &
\raisebox{-\height}{\includegraphics{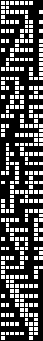}} &
\raisebox{-\height}{\includegraphics{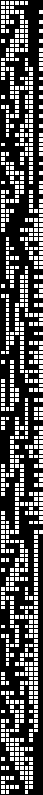}} &
\raisebox{-\height}{\includegraphics{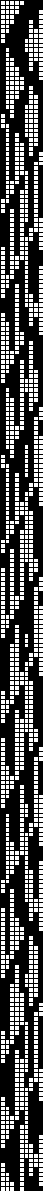}} &
\raisebox{-\height}{\includegraphics{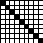}} &
\raisebox{-\height}{\includegraphics{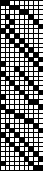}} &
\raisebox{-\height}{\includegraphics{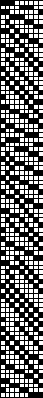}} &
\raisebox{-\height}{\includegraphics{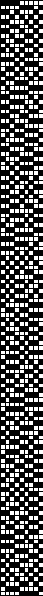}}
\end{tabular}
}
\caption{Hamilton cycles in (a1)--(a4)~bipartite Kneser graphs $H_{9,k}$, for $k=1,2,3,4$;
% Constructed as in the M\"utze/Su paper, with suitable middle levels cycles for the induction basis (these are computed on combos.org, plus swapping of two bits so that the cycles visit a certain triple of vertices consecutively).
and (b1)--(b4)~Kneser graphs $K_{9,k}$, for $k=1,2,3,4$.
% (b1) constructed as in the Nummenpalo/M\"utze/Walczak paper
% (b2)-(b4) constructed as in the Merino/M\"utze/Namrata paper
Vertices are in their bitstring representation ($0=\text{white}$ and $1=\text{black}$).
}
\vspace{-24mm}
\label{fig:kneser}
\end{wrapfigure}
The method of M\"utze and Su extends Havel's idea as follows:
In addition to the cycle~$C_{n,k}$ in~$Q_{n,k}$, we maintain a set of vertex-disjoint paths~$P_{n,k}$ in~$Q_n$, shown in blue in Figure~\ref{fig:bkneser}, each of which starts at a vertex of~$C_{n,k}$ in level~$k+1$ and ends at a vertex in level~$n-k$, and visits only one vertex of each level in between.
The number of such paths equals the number of vertices in levels~$k$ or~$n-k$, namely $\binom{n}{k}=\binom{n}{n-k}$, i.e., these paths visit all vertices in level~$n-k$, but they skip vertices in the levels below that.
For $n=2k+1$ the paths~$P_{n,k}$ have no edges, so this strengthening is again the middle levels conjecture.
For $n>2k+1$ the cycles~$C_{n,k}$ and paths~$P_{n,k}$ can be constructed following a very similar inductive approach as before, by partitioning vertices according to the last bit and gluing together two copies of the structures obtained by induction.

From the cycle~$C_{n,k}$ and the paths~$P_{n,k}$ we construct a Hamilton cycle in~$H_{n,k}$ as follows:
We replace each of the vertices of the cycle~$C_{n,k}$ in level~$k+1$, which is the starting vertex of some path from~$P_{n,k}$, by the other end vertex of this path in level~$n-k$.
This gives a cyclic sequence in which all vertices in level~$k$ are interleaved with all vertices in level~$n-k$, with the additional property that the predecessor and successor of any level-$k$ vertex are reachable from it by a path in~$Q_n$ that moves up from level~$k$ to level~$n-k$.
As moving up along a path in~$Q_n$ corresponds to moving to a superset, this sequence is indeed a Hamilton cycle in~$H_{n,k}$.
This cycle has the remarkable additional closeness property that any two consecutive $k$-sets differ only in the exchange of a single element (as they have a common neighbor in level~$k+1$).
This construction is illustrated in Figure~\ref{fig:kneser}~(a1)--(a4) for the cases~$n=9$ and $k=1,\ldots,4$ ($k=4$ is a solution to the middle levels conjecture).

\section{Kneser graphs}

For integers~$k\geq 1$ and~$n\geq 2k+1$, the \emph{Kneser graph~$K_{n,k}$} has as vertices all $k$-element subset of~$[n]$, and edges between any two disjoint sets.
Kneser graphs have many interesting properties, for example, their chromatic number was shown to be $n-2k+2$ by Lov\'asz using topological methods, and their independence number is $\binom{n-1}{k-1}$ by the famous Erd{\H o}s-Ko-Rado theorem.
Kneser graphs are clearly vertex-transitive, and $H_{n,k}$ is the \emph{bipartite double cover} of~$K_{n,k}$, i.e., we take two copies of~$K_{n,k}$ and replace every corresponding pair of edges inside the copies by the two `diagonal' cross edges between them.
Consequently, if $K_{n,k}$ admits a Hamilton cycle, then~$H_{n,k}$ admits a Hamilton path or cycle.
Indeed, given a Hamilton cycle~$C=(X_1,\ldots,X_\ell)$ in~$K_{n,k}$, where $\ell=\binom{n}{k}$, we define $\ol{X_i}:=[n]\setminus X_i$, and we consider the two sequences~$P:=(X_1,\ol{X_2},X_3,\ol{X_4},\ldots)$ and $P':=(\ol{X_1},X_2,\ol{X_3},X_4,\ldots)$, both of length~$\ell$, obtained by complementing every even- or odd-indexed set in~$C$, respectively.
The last entries of~$P$ and~$P'$ are $X_\ell$ and~$\ol{X_\ell}$, respectively, if $\ell$ is odd, and vice versa if $\ell$ is even.
Furthermore, as $X_i\cap X_{i+1}=\emptyset$ we have $X_i\seq \ol{X_{i+1}}$ and $\ol{X_i}\supseteq X_{i+1}$, so $P$ and~$P'$ are paths in~$H_{n,k}$.
Furthermore, if $\ell$ is odd, then their concatenation~$PP'$ is a Hamilton cycle in~$H_{n,k}$.
On the other hand, if $\ell$ is even, then the two end vertices of~$P$ are adjacent and the two end vertices of~$P'$ are adjacent, so these two disjoint cycles in~$H_{n,k}$ can be joined to a Hamilton path.

The sparsest Kneser graphs are obtained when $n=2k+1$, and they are also known as \emph{odd graphs~$O_k:=K_{2k+1,k}$}.
The odd graph~$O_2=K_{5,2}$ is the Petersen graph shown in Figure~\ref{fig:petersen}, which does not have a Hamilton cycle, but only a Hamilton path.
The graph~$O_3=K_{7,3}$ is shown in Figure~\ref{fig:K73}.
The conjecture that~$O_k$ for $k\geq 3$ has a Hamilton cycle was raised in the 1970s, even before the middle levels conjecture, in papers by Meredith and Lloyd~\cite{MR0321782}, and by Biggs~\cite{MR556008}.
By our earlier observation, the Hamiltonicity of~$K_{n,k}$ implies it for~$H_{n,k}$.
In particular, Hamiltonicity of the odd graphs implies the middle levels conjecture.
Consequently, Kneser graphs attracted a lot of attention, and there was a long line of research on proving that sufficiently dense Kneser graphs~$K_{n,k}$, i.e., those where $n$ is large w.r.t.~$k$, admit a Hamilton cycle.

\subsection{The Chen-F\"uredi construction via Baranyai's partition theorem}

\begin{figure}[b!]
\begin{tabular}{cc}
\raisebox{4.5mm}{\includegraphics{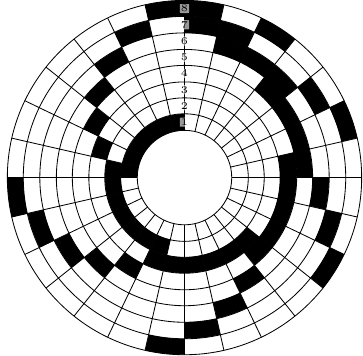}} &
\includegraphics{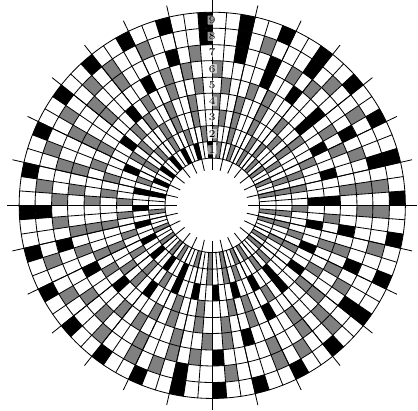} \\
(a) & (b)
\end{tabular}
\caption{(a) Exchange Gray code for 2-element subsets of~[8] obtained by restricting the binary reflected Gray code in~$Q_8$ to level~2; (b) Hamilton cycle in~$K_{9,3}$ via the Chen-F\"uredi construction.
Each indicated triple of vertices is a partition of~$[9]$, the first two are colored gray, and the last one is colored black and corresponds to the set from~(a) by adding the last element (extra outermost black bit).}
\label{fig:baranyai}
\end{figure}

In~2002, Chen and F\"uredi~\cite{MR1883565} found a particularly nice proof that~$K_{n,k}$ has a Hamilton cycle when $n=pk$ for some integer~$p\geq 3$.
The first ingredient of their proof is \emph{Baranyai's partition theorem}, which states that the $\binom{n}{k}$ vertices of the Kneser graph can be partitioned into $\binom{n}{k}/p$ groups of size~$p$ such that the vertices in each group are a partition of~$[n]$, i.e., they are pairwise disjoint and together cover~$[n]$.
% For $k=2$ this is saying that the edges of the complete graph on $n=2p$ vertices can be partitioned into edge-disjoint perfect matchings, which is Walecki's theorem.
The second ingredient is a method to list all $k$-element subset of~$[n]$ in such a way that any two consecutive sets~$X,Y$ differ in an element exchange, i.e., $Y=(X\setminus\{i\})\cup\{j\}$ for some $i,j\in[n]$.
It is well known that such a listing can be obtained from the binary reflected Gray code for~$Q_n$ by restricting it to the vertices in level~$k$, i.e., we simply delete from the full listing all vertices not in level~$k$; see Figure~\ref{fig:baranyai}~(a).

To prove that $K_{n,k}$ with $n=pk$ and $p\geq 3$ has a Hamilton cycle, Chen and F\"uredi first apply Baranyai's theorem, which partitions all vertices of~$K_{n,k}$ into $\ell:=\binom{n}{k}/p$ groups of size~$p$, such that each group is a partition of~$[n]$.
Let $X_1^i,\ldots,X_p^i$ be the sets in the $i$th group, for every $i\in[\ell]$.
Each of those groups forms a clique in the Kneser graph, i.e., it can be traversed in any order.
Without loss of generality we may assume that $n\in X_p^i$, i.e., the element $n$ is contained in the last set of each group.
Furthermore, by the aforementioned Gray code result, we can assume that $X_p^1,X_p^2,\ldots,X_p^\ell$ are ordered so that any two consecutive sets differ in an element exchange, i.e., $X_p^1\setminus\{n\},\ldots,X_p^\ell\setminus\{n\}$ forms an exchange Gray code for all $(k-1)$-element subsets of~$[n-1]$.
Let $x^i\in[n]$ be the element in~$X_p^i$ that is not contained in~$X_p^{i+1}$ (the indices~$i$ are considered modulo~$\ell$).
We may also assume that $X_1^{i+1}$ does not contain~$x^i$, otherwise the sets $X_1^{i+1},\ldots,X_{p-1}^{i+1}$ can be reordered appropriately, which is possible as there are at least $p-1\geq 2$ of them.
It follows that $X_p^i\cap X_1^{i+1}=\emptyset$ and consequently, $X_1^1,\ldots,X_p^1,\,X_1^2,\ldots,X_p^2,\,\ldots,\,X_1^\ell,\ldots,X_p^\ell$ is the desired Hamilton cycle in~$K_{n,k}$; see Figure~\ref{fig:baranyai}~(b).
The method via Baranyai partitions was later refined by Chen~\cite{MR1778200} to establish Hamiltonicity of all~$K_{n,k}$ with $n\geq 2.62k+1$.

\subsection{Settling the odd graphs via a Chung-Feller bijection}

\begin{figure}
\includegraphics{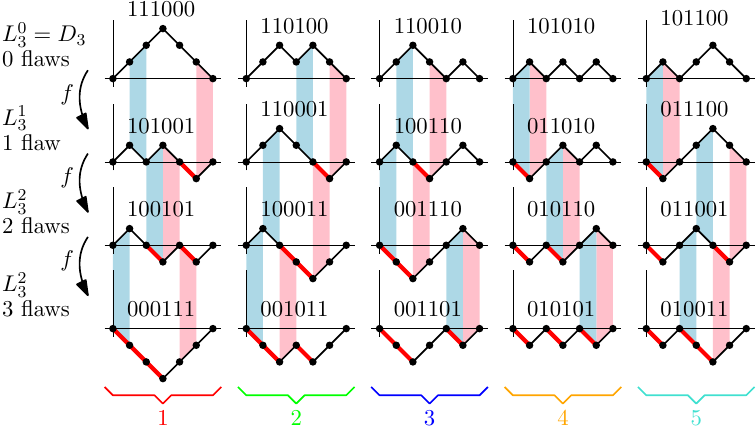}
\caption{Proof of the Chung-Feller theorem for $k=3$ via the minimum change bijection~$f$.
The two transposed bits are highlighted by vertical bars.
Each column produces one of the five cycles of the factor in~$O_3$ shown in Figure~\ref{fig:K73}.}
\label{fig:chung}
\end{figure}

\begin{figure}
\includegraphics{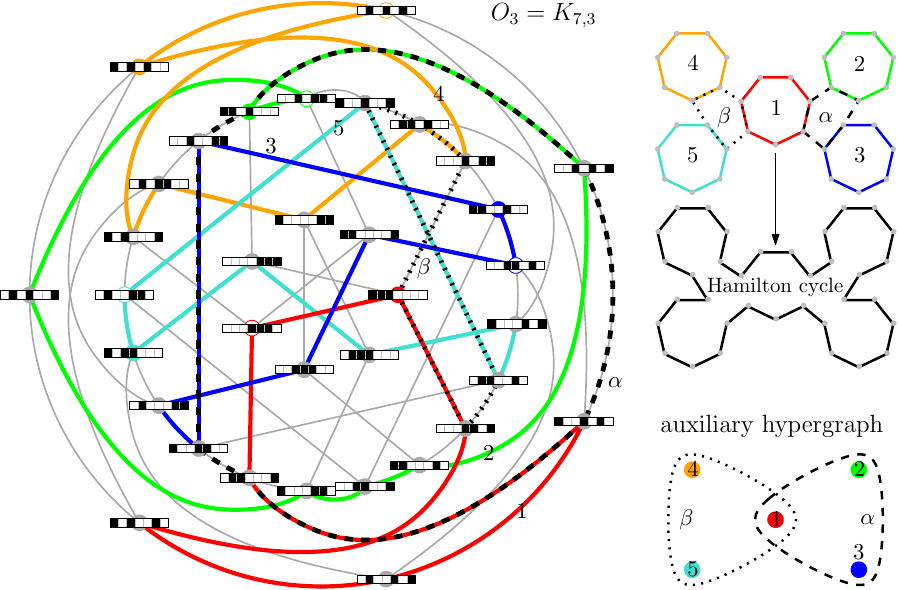}
\caption{Cycle factor in~$O_3$ obtained from the bijection~$f$ in Figure~\ref{fig:chung}, and corresponding gluing cycles to turn it into a Hamilton cycle.}
\label{fig:K73}
\end{figure}

In 2021, M\"utze, Nummenpalo, and Walczak~\cite{MR4273468} proved that the sparsest Kneser graphs, namely the odd graphs~$O_k$ for all~$k\geq 3$ have a Hamilton cycle.
The starting point of their proof is the well-known \emph{Chung-Feller theorem}.
A \emph{flaw} in a bitstring~$x$ is a prefix of~$x$ ending with~0 that has strictly less 1s than~0s; flaws are drawn as red steps in Figure~\ref{fig:chung}.
We write $L_k$ for the unique middle level~$k$ of~$Q_{2k}$, i.e., all bitstrings of length~$2k$ with exactly $k$ many~1s, and we partition~$L_k$ into sets~$L_k^e$ for $e=0,\ldots,k$ according to the number~$e$ of flaws.
In particular, $L_k^0=D_k$ are Dyck words, and $L_k^k$ are complemented Dyck words.
The Chung-Feller theorem asserts that $|L_k^0|=|L_k^1|=\cdots=|L_k^k|=\frac{1}{k+1}\binom{2k}{k}=C_k$, i.e., the number of strings is the same independently of the number of flaws, and it is the $k$th Catalan number.
It is not hard to prove this by establishing a bijection $f:L_k^e\rightarrow L_k^{e+1}$.
In 2018, M\"utze, Standke and Wiechert presented a new proof, using a bijection~$f$ that has the following additional properties; see Figure~\ref{fig:chung}:
$f$ only transposes two bits (0 and~1), for any $x\in L_k$ we have that $f^k(x)$ is the complement of~$x$, i.e., every bit is transposed (and thus complemented) exactly once when applying the bijection $k$~times, and the unique neighbors $\wh{x}:=x\cup f(x)$ of~$x$ and~$f(x)$ in level~$k+1$ of~$Q_{2k}$ for $x\in L_k$ are all distinct and together cover precisely this level.
We can thus build vertex-disjoint paths $(x,\wh{x},f(x),\wh{f(x)},f^2(x),\wh{f^2(x)},\ldots,f^k(x))$, all of length~$2k$, that together cover~$Q_{2k,k}$ and that connect pairs of Dyck words and their complements (the length of each path is the number of its edges, which is one less than the number of vertices).
Appending a 0-bit to all vertices and taking complements of the resulting vertices in level~$k+1$ yields a cycle factor in~$O_k$ which has $C_k$ many cycles of the same length~$2k+1$; see Figure~\ref{fig:K73}.

The proof of M\"utze, Nummenpalo, and Walczak is completed by gluing the cycles of this factor together via 6-cycles and 8-cycles, which glue together 3 or 4 cycles, respectively, from the factor at a time.
The main technical difficulty is that the corresponding auxiliary graph is now a hypergraph with hyperedges of cardinality~3 or~4, respectively; see the bottom right of Figure~\ref{fig:K73} (only 3-hyperedges are present in the figure).
To obtain a Hamilton cycle in~$O_k$, we seek a so-called \emph{loose} spanning tree in the auxiliary hypergraph, i.e., a spanning tree in which any two hyperedges overlap at most in a singleton.
For this it is not enough to prove that the hypergraph is connected, as there are connected hypergraphs that do not admit any loose spanning tree.
Instead the paper constructs one particular loose spanning tree.
The resulting Hamilton cycle in~$O_4$ is shown in Figure~\ref{fig:kneser}~(b4).

\subsection{Johnson's inductive construction}

In~2011, Johnson~\cite{MR2836824} devised an inductive construction for Hamilton cycles in Kneser graphs.
Specifically, he showed that $K_{n,k}$ for $n=2k+s$ with even~$s$ has a Hamilton cycle provided that the smaller Kneser graphs $K_{2\ell+s/2,\ell}$ have a Hamilton cycle for all $1\leq\ell\leq \lfloor k/2\rfloor$ (or they are the Petersen graph~$K_{5,2}$).

His construction works for a ground set of even cardinality~$n=2k+s$ by partitioning it into fixed pairs~$\{2i-1,2i\}$ for $i=1,\ldots,n/2$, and by considering the possible intersection patterns of subsets with those pairs.
Specifically, we associate a set $X\seq[n]$ with a tuple~$(X(1),\ldots,X(n/2))\in \{-1,0,1,2\}^{n/2}$ by defining $X(i):=0,-1,1$, or~$2$ if $X\cap \{2i-1,2i\}$ equals $\emptyset,\{2i-1\},\{2i\}$, or $\{2i-1,2i\}$, respectively.
In the simplest case, the subsets~$X$ intersect the pairs in either 0 or 2 elements, i.e., $X(i)\in\{0,2\}$ for all $i\in[n/2]$.
In this case, for $\ell=k/2$ a Hamilton cycle in~$K_{k+s/2,k/2}$ can be lifted to a cycle in~$K_{n,k}=K_{2k+s,k}$ by replacing every element~$i$ by the pair~$\{2i-1,2i\}$.
More generally, consider subsets~$X$ which intersect all but a fixed set of $t$ pairs in~0 or~2 elements, i.e., $X(i)\in\{-1,1\}$ for a fixed $t$-set of indices~$i\in[n/2]$.
An edge in $K_{k+s/2-t,k/2-t/2}$ together with a $t$-set in~$[n/2]$ lifts to a set of edges involving subsets which intersect all but this fixed set of $t$ pairs in~0 or~2 elements.
For example, we have many edges (one for each pattern of $\pm1$s) of the form $((0,0,2,0,2,2,0,\pm1,\pm1), (2,2,0,0,0,0,2,\mp1,\mp1))$ in~$K_{18,8}$ that arise from the edge $(\{3,5,6\},\{1,2,7\})$ in $K_{7,3}$ with $\{8,9\}$ as the special $t$-set.
Using this idea, a Hamilton cycle in $K_{k+s/2-t,k/2-t/2}$ lifts to a cycle (possibly of double the length) consisting of sets of this type in~$K_{n,k}$.
With some care, one can join together these cycles corresponding to different $t$-sets, and then for different values of~$t$ to give a Hamilton cycle in~$K_{n,k}$.

Combining Johnson's result with the solution for the sparsest case $n=2k+1$ presented in the previous section, we obtain that $K_{2k+2^a,k}$ has a Hamilton cycle for all $k\geq 3$ and $a\geq 0$.
This settles in particular the second-sparsest case~$K_{2k+2,k}$.

\subsection{Settling the remaining cases via Greene-Kleitman parenthesis matching and gliders}

\begin{figure}
\includegraphics[page=1]{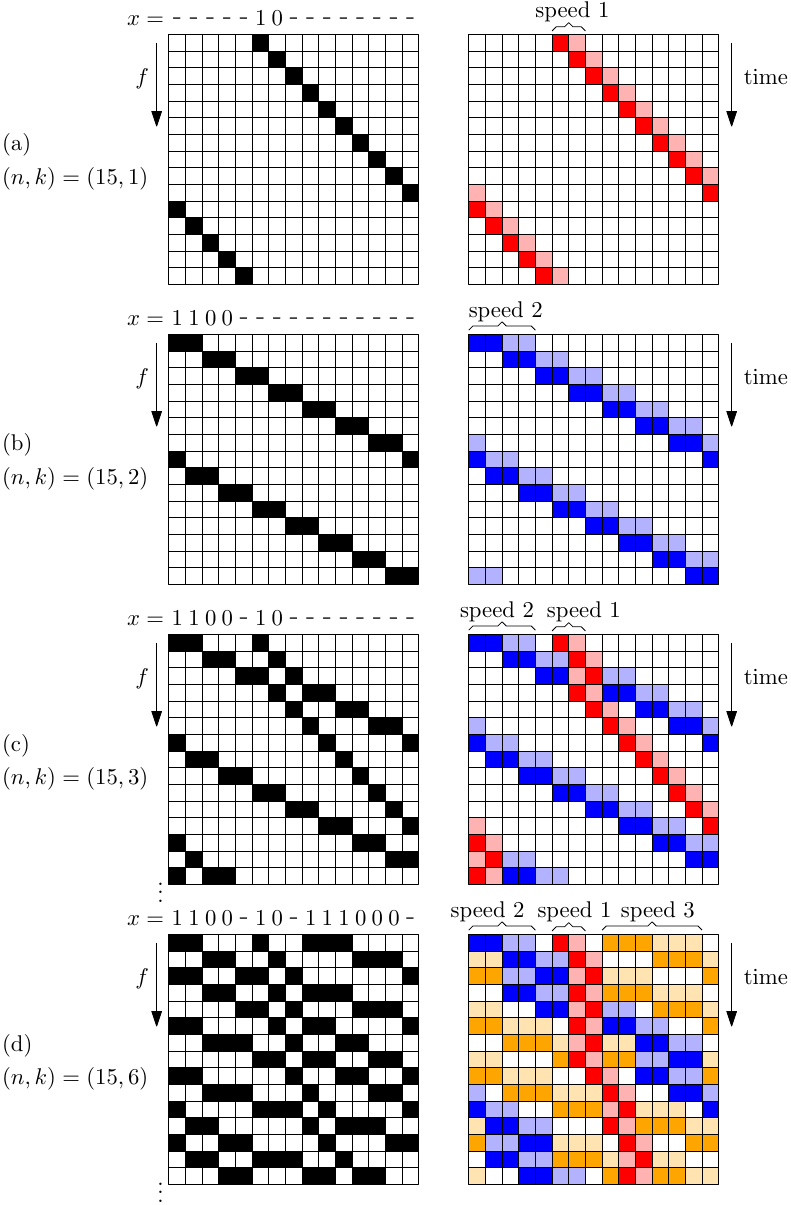}
\caption{Cycles in different Kneser graphs~$K_{n,k}$ constructed by parenthesis matching.
The cycles in~(a) and~(b) are shown completely, whereas in~(c) and~(d) only the first 15 vertices are shown.
When applying parenthesis matching to~$x$, unmatched 0s are printed as~$\hyph$.
The right hand side shows the interpretation of certain groups of bits as gliders, and their movement over time.
Matched bits belonging to the same glider are colored in the same color, with the opaque filling given to 1-bits, and the transparent filling given to 0-bits.
(a) one glider of speed~1; (b) one glider of speed~2; (c) two gliders with speeds~1 and~2 that participate in an overtaking; (d) three gliders of speeds~1, 2 and~3 that participate in multiple overtakings.
Animations of these examples are available at~\cite{gliders}.
}
\label{fig:gliders}
\end{figure}

In 2022, Merino, M\"utze, and Namrata~\cite{merino-muetze-namrata:22} proved that $K_{n,k}$ for $n\geq 2k+3$ has a Hamilton cycle, which combined with the results from the previous two sections completely settles the problem for Kneser graphs.
Their proof starts with a new cycle factor in~$K_{n,k}$, which is constructed using the following simple rule based on \emph{parenthesis matching}, a technique that was pioneered by Greene and Kleitman~\cite{MR0389608} in the context of symmetric chain partitions of the Boolean lattice:
We consider vertices of~$K_{n,k}$ as bitstrings, and we interpret the 1s in~$x$ as opening brackets and the 0s as closing brackets, and we match closest pairs of opening and closing brackets in the natural way, which will leave some 0s unmatched.
This matching is done \emph{cyclically} across the boundary of~$x$, i.e., $x$ is considered as a cyclic string.
We write $f(x)$ for the vertex obtained from~$x$ by complementing all matched bits, leaving the unmatched bits unchanged.
Note that $x$ and~$f(x)$ have no 1s at the same positions, implying that $(x,f(x))$ is an edge in the Kneser graph.
Furthermore, $f$ is invertible and and~$f^2(x)\neq x$, so the union of all edges~$(x,f(x))$ is a collection of disjoint cycles that together visit all vertices of~$K_{n,k}$; see Figure~\ref{fig:gliders}.

The next step is to understand the structure of the cycles generated by~$f$.
Interestingly, the evolution of a bitstring~$x$ under repeated applications of~$f$ can be described by a kinetic system of multiple gliders that move at different speeds and that interact over time, somewhat reminiscent of the gliders in Conway's Game of Life.
Specifically, each application of~$f$ is viewed as one unit of time moving forward.
Furthermore, we partition the matched bits of~$x$ into groups, and each of these groups is called a \emph{glider}.
A glider has a \emph{speed} associated to it, which is given by the number of 1s in its group.
For example, in the cycle shown in Figure~\ref{fig:gliders}~(a), there is a single matched~1 and the corresponding matched~0, and together these two bits form a glider of speed~1 that moves one step to the right in every time step.
Applying $f$ means going down to the next row in the picture, so the time axis points downwards.
Similarly, in Figure~\ref{fig:gliders}~(b), there are two matched~1s and the corresponding two matched~0s, and together these four bits form a glider of speed~2 that moves two steps to the right in every time step.
As we see from these examples, a single glider of speed~$v$ simply moves uniformly, following the basic physics law
\begin{equation*}
s(t)=s(0)+v\cdot t,
\end{equation*}
where $t$ is the time (i.e., the number of applications of~$f$) and $s(t)$ is the position of the glider in the bitstring as a function of time.
The position~$s(t)$ has to be considered modulo~$n$, as bitstrings are considered as cyclic strings and the gliders hence wrap around the boundary.
The situation gets more interesting and complicated when gliders of different speeds interact with each other.
For example, in Figure~\ref{fig:gliders}~(c), there is one glider of speed~2 and one glider of speed~1.
As long as these groups of bits are separated, each glider moves uniformly as before.
However, when the speed~2 glider catches up with the speed~1 glider, an overtaking occurs.
During an overtaking, the faster glider receives a boost, whereas the slower glider is delayed.
This can be captured by augmenting the corresponding equations of motion by introducing additional terms, making them non-uniform.
In the simplest case of two gliders of different speeds, the equations become
\begin{align*}
s_1(t)&=s_1(0)+v_1\cdot t-2v_1 c_{1,2}, \\
s_2(t)&=s_2(0)+v_2\cdot t+2v_1 c_{1,2},
\end{align*}
where the subscript~1 stands for the slower glider and the subscript~2 stands for the faster glider, and the additional variable~$c_{1,2}$ counts the number of overtakings.
Note that the terms~$2v_1c_{1,2}$ occur with opposite signs in both equations, capturing the fact that the faster glider is boosted by the same amount that the slower glider is delayed.
This can be seen as `energy conservation' in the system of gliders.
For more than two gliders, the equations of motion can be generalized accordingly, by introducing additional overtaking counters between any pair of gliders.
From those equations of motions, important properties of the cycles can be extracted via combinatorial and algebraic arguments.
One such property is that the number of gliders and their speeds are invariant along each cycle.
For example, in Figure~\ref{fig:gliders}~(d), every bitstring along this cycle has three gliders of speeds~1, 2 and~3.
For the reader's entertainment, we programmed an interactive animation of gliders over time, and we encourage experimentation with this code, which can be found at~\cite{gliders}.

The last step of the proof joins the cycles of this factor via gluing 4-cycles (this is where the assumption $n\geq 2k+3$ is used).
Specifically, the gluing cycles join pairs of cycles whose sets of glider speeds differ in a small modification, changing the speed of one glider by~$-1$ and the speed of another by~$+1$.
We thus obtain a combinatorial interpretation of the gluings.
To prove that all cycles of the factor can be joined to a single Hamilton cycle, it is argued that all cycles can be joined to one particular cycle in the factor, by considering the speed sets of gliders as number partitions, and by arguing that these partitions increase lexicographically along suitable gluings.
The Hamilton cycles in~$K_{9,k}$ for $k=1,2,3$ resulting from this proof are shown in Figure~\ref{fig:kneser}~(b1)--(b3).

\section{Generalized Johnson graphs}

The \emph{generalized Johnson graph $J_{n,k,s}$} has as vertices all $k$-element subsets of~$[n]$, and an edge between any two sets whose intersection has size exactly~$s$.
It is defined for integers~$k\geq 1$, $0\leq s<k$ and $n\geq 2k-s+\mathbf{1}_{s=0}$, where $\mathbf{1}_{s=0}$ denotes the indicator function that equals~1 if $s=0$ and~0 otherwise.
For $s=0$ we obtain Kneser graphs ($s=0$), and for $s=k-1$ we obtain Johnson graphs as special cases.
By taking complements, we see that $J_{n,k,s}$ is isomorphic to~$J_{n,n-k,n-2k+s}$.
Chen and Lih~\cite{MR888679} conjectured in~1987 that all graphs~$J_{n,k,s}$ admit a Hamilton cycle except the Petersen graph~$J_{5,2,0}=J_{5,3,1}$.
This includes Hamiltonicity of the corresponding bipartite double covers, in particular a solution to the middle levels conjecture, which was the starting point of this article.

In fact, already Chen and Lih observed that~$J_{n,k,s}$ can be partitioned into two subgraphs isomorphic to~$J_{n-1,k,s}$ and~$J_{n-1,k-1,s-1}$ (split vertices according to containment of some fixed element, $n$ say), so if these two graphs have a Hamilton cycle, then we can glue them via a 4-cycle and obtain a Hamilton cycle in~$J_{n,k,s}$.
To complete the proof, it remains to observe that if~$J_{n,k,s}$ is a generalized Johnson graph, then either it is a Kneser graph, or $J_{n-1,k,s}$ and $J_{n-1,k-1,s-1}$ are both generalized Johnson graphs.
Using the results for Kneser graphs from the previous section we thus obtain Hamiltonicity for all generalized Johnson graphs by induction.

There is another closely related and heavily studied class of vertex-transitive graphs called \emph{generalized Kneser graphs~$K_{n,k,s}$}.
This graph has as vertices all $k$-element subsets of~$[n]$, and an edge between any two sets whose intersection has size \emph{at most}~$s$.
Clearly, $J_{n,k,s}$ is a spanning subgraph of~$K_{n,k,s}$, so the Hamiltonicity of~$K_{n,k,s}$ is immediate.

\section{What's next?}

Generalized Johnson graphs are the most general family of graphs defined by intersecting set systems and thus a natural end to our story.
With regards to Lov\'asz' conjecture, many other families of vertex-transitive graphs await to be tested for Hamiltonicity, in particular Cayley graphs.
We may also ask how many distinct Hamilton cycles a vertex-transitive graph admits.
In fact, the proofs via gluings of cycle factors for the middle levels graph~$M_k$ and for the odd graph~$O_k$ discussed in this article yield double-exponentially (in $k$) many distinct Hamilton cycles, and the trivial upper bound of $n!$ for the number of Hamilton cycles of an $n$-vertex graphs also yields a double-exponential function in~$k$.
A much harder problem is to find multiple \emph{edge-disjoint} Hamilton cycles.
Biggs~\cite{MR556008} conjectured that the odd graph~$O_k$ can be partitioned into $\lfloor (k+1)/2\rfloor$ edge-disjoint Hamilton cycle for~$k\geq 3$.
We do not even know two edge-disjoint Hamilton cycles in the middle levels graph~$M_k$.
Katona conjectured that the Kneser graph~$K_{n,k}$ contains the $r$th power of a Hamilton cycle, where $r:=\lfloor n/k\rfloor-2$, and this may even be true for $r:=\lceil n/k\rceil-2$. 
Also, for the graphs considered in this article, we may ask whether they are \emph{Hamilton-connected}, i.e., they admit a Hamilton path between any two prescribed end vertices.
For bipartite graphs, we may ask whether they are \emph{Hamilton-laceable}, i.e., they admit a Hamilton path between any two prescribed end vertices, one from each partition class.
In fact, the middle levels graph~$M_k$ for $k\geq 2$ was shown to be Hamilton-laceable.
Another generalization is to consider the \emph{cycle spectrum}, which is the set of all possible cycle lengths in a graph.
For example, does the middle levels graph~$M_k$ admit cycles of all possible even lengths starting from~6 up to the number of vertices?

From an algorithmic point of view, one may ask which of the cycles described in this article can be computed efficiently?
A satisfactory answer to this question is only known for the middle levels graph~$M_k$, while all the other known constructions present fundamental obstacles to such algorithms.
Furthermore, what about simple descriptions of Hamilton cycles, similar in flavor to Williams' greedy description of the binary reflected Gray code mentioned in Section~\ref{sec:cube}?
Even the simplest known solution of the middle levels conjecture is much more complicated than this.

There are many other intriguing problems about the interaction of different structures, such as matchings and cycles, in vertex-transitive graphs.
For example, Ruskey and Savage asked whether every matching in the hypercube can be extended to a Hamilton cycle.
For the case of perfect matchings this was answered affirmatively by Fink~\cite{MR2354719}.
Also, Kotzig's question on perfect 1-factorizations of the complete graph comes to mind naturally.
A \emph{perfect 1-factorization} is a decomposition of the edge set of a graph into perfect matchings, such that the union of any two of them forms a Hamilton cycle.

\bibliographystyle{alpha}
\bibliography{refs}

\end{document}